\documentclass[a4j,11pt]{article}
\usepackage{amsmath,amsthm,amssymb,amscd,ascmac}
\usepackage{graphicx, mathrsfs}
\theoremstyle{definition}
\numberwithin{equation}{section}
\newtheorem{thm}{Theorem}[section]

\newtheorem{lem}[thm]{Lemma}
\newtheorem{cor}[thm]{Corollary}

\newcommand{\tr}{\mathop{\mathrm{tr}}\,}
\renewcommand{\det}{\mathop{\mathrm{det}}\,}

\theoremstyle{definition}

\newtheorem{rem}[thm]{Remark}

\usepackage{bm}

\begin{document}
\title{An explicit formula of powers of the $2\times 2$ quantum matrices and its applications}
\author{Genki Shibukawa\thanks{
Dedicated to T. Umeda and M. Wakayama for their 66th birthdays.}}
\date{
\small MSC classes\,:\,16T20, 33C45, 81R50}
\pagestyle{plain}

\maketitle

\begin{abstract}
We present an explicit formula of the powers for the $2\times 2$ quantum matrices, that is a natural quantum analogue of the powers of the usual $2\times 2$ matrices. 
As applications, we give some non-commutative relations of the entries of the powers for the $2\times 2$ quantum matrices, which is a simple proof of the results of Vokos-Zumino-Wess (1990).
\end{abstract}

\section{Introduction}
Let $A$ be a $2\times 2$ matrix over a fixed base field $k$ and $a,b,c,d \in k$ are entries of $A$:
$$
A=
   \begin{pmatrix}
   a & b \\
   c & d
   \end{pmatrix}.
$$
For any positive number $n$, the following explicit formula of $A^{n}$ holds:
\begin{align}
\label{eq:powers origin}
A^{n}
   =\begin{pmatrix}
   af_{n-1}\left(\tr{A} \right)-\det{A} f_{n-2}\left(\tr{A} \right) & bf_{n-1}\left(\tr{A} \right) \\
   cf_{n-1}\left(\tr{A} \right) & df_{n-1}\left(\tr{A} \right)-\det{A} f_{n-2}\left(\tr{A} \right) 
  \end{pmatrix}.
\end{align}
Here $f_{n}(x)$ is the polynomial of degree $n$ defined by
\begin{align}
f_{n}(x)
   =
   f_{n}(x,y)
   :=
   \sum_{l=0}^{\left\lfloor \frac{n}{2}\right\rfloor}
   (-1)^{l}\binom{n-l}{l}x^{n-2l}y^{l}, \quad 
f_{-1}(x):=0, 
\end{align}
which is the Chebyshev polynomial of the second kind:
$$
U_{n}(x):=\sum_{l=0}^{\left\lfloor \frac{n}{2}\right\rfloor}(-1)^{l}\binom{n-l}{l}x^{n-2l}.
$$
In fact this formula is well-known in linear algebra, and its proof is easy by induction and the recurrence relation for $f_{n}(x)$:
\begin{align}
\label{eq:rec of f}
f_{n+1}(x)=xf_{n}(x)-yf_{n-1}(x).
\end{align}

In this note, we give a quantum analogue of the formula (\ref{eq:powers origin}) and its applications. 
First, we review some fundamental objects and facts on quantum matrix or groups \cite{M}, \cite{T} relevant to the main results of this paper.

We call $A$ a $2\times 2$ ($q$-)quantum matrix if its entries satisfy the following relations:
\begin{align}
& ab=qba, \quad ac=qca, \quad ad-da=(q-q^{-1})bc, \nonumber \\
\label{eq:quantum rel} 
& bc=cb, \quad bd=qdb, \quad cd=qdc, \tag{$R_{q}$}
\end{align}
Here $q$ is a central indeterminate. 
A quantum analogue of the coordinate ring $\mathcal{A}_{q}(Mat(2))$ is the algebra generated by $a,b,c,d$ and $q$ which is a typical example of quantum groups.

The quantum adjoint matrix of any quantum matrix $A$ is defined as
\begin{align}
\hat{A}
   =\begin{pmatrix}
   \hat{a} & \hat{b} \\
   \hat{c} & \hat{d}
   \end{pmatrix}
   :=\begin{pmatrix}
   d & -q^{-1}b \\
   -qc & a
   \end{pmatrix}.
\end{align}
By the definition of $\hat{A}$ and the relations (\ref{eq:quantum rel}), the quantum adjoint matrix $\hat{A}$ satisfies the relations $(R_{q^{-1}})$:
\begin{align*}
& \hat{a}\hat{b}=q^{-1}\hat{b}\hat{a}, \quad \hat{a}\hat{c}=q^{-1}\hat{c}\hat{a}, \quad \hat{a}\hat{d}-\hat{d}\hat{a}=(q^{-1}-q)\hat{b}\hat{c}, \\
& \hat{b}\hat{c}=\hat{c}\hat{b}, \quad \hat{b}\hat{d}=q^{-1}\hat{d}\hat{b}, \quad \hat{c}\hat{d}=q^{-1}\hat{d}\hat{c}.
\end{align*}
Hence, the relations (\ref{eq:quantum rel}) are equivalent to 
\begin{align}
\label{eq:AhatA}
bc=cb, \quad A\hat{A}=\hat{A}A=
\begin{pmatrix}
\delta  & 0 \\
0 & \delta 
\end{pmatrix}
=\delta E_{2},
\end{align}
where $E_{2}$ is the $2\times 2$ identity matrix, and $\delta :=ad-qbc=da-q^{-1}bc$ is the quantum determinant of $A$ which is a central element of $\mathcal{A}_{q}(Mat(2))$.

For convenience, we introduce a $2\times 2$ matrix 
$$
C
   :=
   \begin{pmatrix}
   q^{\frac{1}{2}} & 0 \\
   0 & q^{-\frac{1}{2}}
   \end{pmatrix}   
$$
and put 
$$
\tau :=\tr{(AC)}=q^{\frac{1}{2}}a+q^{-\frac{1}{2}}d, \quad 
\tau ^{\prime}:=\tr{(C^{-1}A)}=q^{-\frac{1}{2}}a+q^{\frac{1}{2}}d.
$$
Our main results are following.
\begin{thm}
\label{thm:main theorem1}
For any positive integer $n$, we have
\begin{align}
\label{eq:power CH1}
A^{n}
   &=
   AC^{-n+1}f_{n-1}(\tau )
   -C^{-n}\delta f_{n-2}(\tau ) \\
\label{eq:power CH2}
   &=
   f_{n-1}(\tau^{\prime})C^{n-1}A
   -f_{n-2}(\tau^{\prime})\delta C^{n},
\end{align}
where 
$$
f_{n}(\tau )
   :=
   f_{n}(\tau, \delta )
   =
   \sum_{l=0}^{\left\lfloor \frac{n}{2}\right\rfloor}
   (-1)^{l}\binom{n-l}{l}\tau ^{n-2l}\delta ^{l}, \quad 
f_{-1}(\tau ):=0.
$$
\end{thm}
Let us put
$$
A^{n}
   =
   \begin{pmatrix}
   a_{n} & b_{n} \\
   c_{n} & d_{n}
   \end{pmatrix}, \quad 
\hat{A}^{m}
   :=
   \begin{pmatrix}
   \hat{a}_{m} & \hat{b}_{m} \\
   \hat{c}_{m} & \hat{d}_{m}
   \end{pmatrix}.
$$
By comparing the entries of $A^{n}$ and (\ref{eq:power CH1}), (\ref{eq:power CH2}), we obtain the following explicit formulas of the entries of $A^{n}$.
\begin{cor}
\label{thm:main theorem2}
For any positive integer $n$, we have
\begin{align}
\label{eq:explicit a}
a_{n}
   &=
   q^{-\frac{n-1}{2}}a
   f_{n-1}(\tau )
   -q^{-\frac{n}{2}}\delta f_{n-2}(\tau )
   =
   q^{\frac{n-1}{2}}
   f_{n-1}(\tau^{\prime })a
   -q^{\frac{n}{2}}\delta f_{n-2}(\tau^{\prime } ), \\
\label{eq:explicit b}
b_{n}
   &=
   q^{\frac{n-1}{2}}b
   f_{n-1}(\tau )
   =
   q^{\frac{n-1}{2}}
   f_{n-1}(\tau^{\prime} )
   b, \\
\label{eq:explicit c}
c_{n}
   &=
   q^{-\frac{n-1}{2}}c
   f_{n-1}(\tau )
   =
   q^{-\frac{n-1}{2}}
   f_{n-1}(\tau^{\prime} )
   c, \\
\label{eq:explicit d}
d_{n}
   &=
   q^{\frac{n-1}{2}}d
   f_{n-1}(\tau )
   -q^{\frac{n}{2}}\delta 
   f_{n-2}(\tau ) 
   =
   q^{-\frac{n-1}{2}}
   f_{n-1}(\tau^{\prime} )
   d
   -q^{-\frac{n}{2}}\delta f_{n-2}(\tau^{\prime}).
\end{align}
\end{cor}
As applications, we give the following results, in particular Theorem \ref{thm:VZW}. 
\begin{cor}
\label{thm:main theorem2-2}
For any positive integer $m$, we have
\begin{align}
\label{eq:explicit hat a}
\hat{a}_{m}
   &=
   q^{\frac{m-1}{2}}\hat{a}
   f_{m-1}(\tau )
   -q^{\frac{m}{2}}\delta f_{m-2}(\tau ) \nonumber \\
   &=
   q^{-\frac{m-1}{2}}
   f_{m-1}(\tau^{\prime })\hat{a}
   -q^{-\frac{m}{2}}\delta f_{m-2}(\tau^{\prime } )
   =
   d_{m}, \\
\label{eq:explicit hat b}
\hat{b}_{m}
   &=
   q^{-\frac{m-1}{2}}\hat{b}
   f_{m-1}(\tau )
   =
   q^{-\frac{m-1}{2}}
   f_{m-1}(\tau^{\prime} )
   \hat{b}
   =
   -q^{-m}b_{m}, \\
\label{eq:explicit hat c}
\hat{c}_{m}
   &=
   q^{\frac{m-1}{2}}\hat{c}
   f_{m-1}(\tau )
   =
   q^{\frac{m-1}{2}}
   f_{m-1}(\tau^{\prime} )
   \hat{c}
   =
   -q^{m}c_{m}, \\
\label{eq:explicit hat d}
\hat{d}_{m}
   &=
   q^{-\frac{m-1}{2}}\hat{d}
   f_{m-1}(\tau )
   -q^{-\frac{m}{2}}\delta 
   f_{m-2}(\tau ) \nonumber \\
   &=
   q^{-\frac{m-1}{2}}
   f_{m-1}(\tau^{\prime} )
   \hat{d}
   -q^{-\frac{m}{2}}\delta f_{m-2}(\tau^{\prime})
   =
   a_{m}.
\end{align}
\end{cor}
\begin{thm}[\cite{VZW}]
\label{thm:VZW}
For any non-negative integers $m$ and $n$, we have
\begin{align}
\label{eq:quantum rel1}
d_{m}a_{n}-q^{-m}b_{m}c_{n}
   &=
   a_{n}d_{m}-q^{m}b_{n}c_{m}
   =
   \begin{cases}
   \delta ^{m}a_{n-m} & (m< n)\\
   \delta ^{n}d_{m-n} & (m\geq n)
   \end{cases}, \\
\label{eq:quantum rel2}
d_{m}b_{n}-q^{-m}b_{m}d_{n}
   &=
   -q^{-m}a_{n}b_{m}+b_{n}a_{m}
   =
   \begin{cases}
   \delta ^{m}b_{n-m} & (m< n)\\
   -q^{n-m}\delta ^{n}b_{m-n} & (m\geq n)
   \end{cases}, \\
\label{eq:quantum rel3}
-q^{m}c_{m}a_{n}+a_{m}c_{n}
   &=
   c_{n}d_{m}-q^{m}d_{n}c_{m}
   =
   \begin{cases}
   \delta ^{m}c_{n-m} & (m< n)\\
   -q^{n-m}\delta ^{n}c_{m-n} & (m\geq n)
   \end{cases}, \\
\label{eq:quantum rel4}
-q^{m}c_{m}b_{n}+a_{m}d_{n}
   &=
   -q^{-m}c_{n}b_{m}+d_{n}a_{m}
   =
   \begin{cases}
   \delta ^{m}d_{n-m} & (m< n)\\
   \delta ^{n}a_{m-n} & (m\geq n)
   \end{cases}, \\
\label{eq:quantum rel5}
b_{n}c_{m}-q^{n-m}c_{n}b_{m}
   &=
   0.
\end{align}
\end{thm}

\section{Proof of Theorem \ref{thm:main theorem1}}
To prove Theorem \ref{thm:main theorem1}, we need a quantum analogue of Cayley-Hamilton theorem.
\begin{lem}[\cite{UW} Lemma 3]
The following formula holds. 
\begin{align}
\label{eq:qCH}
A^{2}
   =
   AC^{-1}\tau -C^{-2}\delta
   =
   \tau^{\prime}CA -\delta C^{2}.
\end{align}
\end{lem}
 
{\bf{Proof of Theorem \ref{thm:main theorem1}}} 
Since (\ref{eq:power CH1}) and (\ref{eq:power CH2}) can be similarly proved, we only prove (\ref{eq:power CH1}). 
These formulas are proved by induction on $n$.

The $n=1$ case is trivial. 
Assume the case of $n$ holds. 
Hence, from the induction hypothesis we have
\begin{align*}
A^{n+1}
   &=
   A A^{n} \\
   &=
   A\left\{AC^{-n+1}f_{n-1}(\tau )
   -C^{-n}\delta f_{n-2}(\tau )\right\} \\
   &=
   A^{2}C^{-n+1}f_{n-1}(\tau )
   -AC^{-n}\delta f_{n-2}(\tau ).
\end{align*}
By Cayley-Hamilton theorem (\ref{eq:qCH}) and the recursion (\ref{eq:rec of f}), we have 
\begin{align*}
A^{n+1}
   &=
   (AC^{-1}\tau -C^{-2}\delta)C^{-n+1}f_{n-1}(\tau )
   -AC^{-n}\delta f_{n-2}(\tau ) \\
   &=
   AC^{-n}(\tau f_{n-1}(\tau )-\delta f_{n-2}(\tau ))
   -C^{-n-1}\delta f_{n-1}(\tau ) \\
   &=
   AC^{-n}f_{n}(\tau )-C^{-n-1}\delta f_{n-1}(\tau ).
\end{align*}
The formula (\ref{eq:power CH2}) can be proved by the similar argument for $A^{n+1}=A^{n}A$. \qed

Corollary \ref{thm:main theorem2} follows from comparing the entries of $A^{n}$ and (\ref{eq:power CH1}), (\ref{eq:power CH2}) immediately. 

\begin{rem}
{\rm{(1)}} By consider the classical limit $q=1$ in Corollary \ref{thm:main theorem2}, we recover the classical result (\ref{eq:powers origin}).\\
{\rm{(2)}} From the recursion (\ref{eq:rec of f}), we derive other expressions of $a_{n}$ and $d_{n}$:
\begin{align}
\label{eq:explicit a2}
a_{n}
   &=
   q^{-\frac{n}{2}}
   f_{n}(\tau )
   -q^{-\frac{n+1}{2}}d f_{n-1}(\tau )
   =
   q^{\frac{n}{2}}
   f_{n}(\tau^{\prime })
   -q^{\frac{n+1}{2}} f_{n-1}(\tau^{\prime } )d, \\
\label{eq:explicit d2}
d_{n}
   &=
   q^{\frac{n}{2}}
   f_{n}(\tau )
   -q^{\frac{n+1}{2}}a
   f_{n-1}(\tau ) 
   =
   q^{-\frac{n}{2}}
   f_{n}(\tau^{\prime} )
   -q^{-\frac{n+1}{2}}f_{n-1}(\tau^{\prime})a.
\end{align}
Since these expressions (\ref{eq:explicit a2}), (\ref{eq:explicit d2}) (and (\ref{eq:explicit b}), (\ref{eq:explicit c})) hold for $n=0$, Corollary \ref{thm:main theorem2} is also true for the case of $n=0$.\\
{\rm{(3)}} Umeda-Wakayama \cite{UW} considered 
\begin{align*}
\tau _{n}&:=\tr{A^{n}C^{n}}=q^{\frac{n}{2}}a_{n}+q^{-\frac{n}{2}}b_{n}, \\
\tau _{n}^{\prime }&:=\tr{C^{-n}A^{n}}=q^{-\frac{n}{2}}a_{n}+q^{\frac{n}{2}}b_{n},
\end{align*}
and pointed out that $\tau _{n}$ and $\tau _{n}^{\prime }$ satisfy the following Fibonacci type equations:
\begin{align}
\label{eq:rec of tau}
\tau _{n+1}=\tau _{n}\tau -\tau _{n-1}\delta , \quad 
\tau _{n+1}^{\prime }=\tau _{n}^{\prime }\tau^{\prime } -\tau _{n-1}^{\prime }\delta .
\end{align}
These equations (\ref{eq:rec of tau}) are equal to the recursion (\ref{eq:rec of f}) of $f_{n}(\tau )$ exactly. 
Hence by $\tau _{1}=\tau$ and $\tau _{1}^{\prime }=\tau^{\prime }$ we have
\begin{align}
\tau _{n}=f_{n}(\tau ), \quad \tau _{n}^{\prime }=f_{n}(\tau ^{\prime }).
\end{align}
\end{rem}

\section{Applications}
We point out that quantum adjoint matrix $\hat{A}$ is a $q^{-1}$-quantum matrix and 
\begin{align*}
\hat{\tau}&:=q^{-\frac{1}{2}}\hat{a}+q^{\frac{1}{2}}\hat{d}=q^{\frac{1}{2}}a+q^{-\frac{1}{2}}d=\tau \\
\hat{\tau}^{\prime}&:=q^{\frac{1}{2}}\hat{a}+q^{-\frac{1}{2}}\hat{d}=q^{-\frac{1}{2}}a+q^{\frac{1}{2}}d=\tau ^{\prime}, \\ 
\hat{\delta}&:=\hat{a}\hat{d}-q^{-1}\hat{b}\hat{c}=da-q^{-1}bc=\delta.
\end{align*}
Then we prove Corollary \ref{thm:main theorem2-2}.

From (\ref{eq:AhatA}), Corollary \ref{thm:main theorem2} and Corollary \ref{thm:main theorem2-2}, we obtain the proof of Theorem \ref{thm:VZW} which is a simple prood of Vokos-Zumino-Wess \cite{VZW}.\\

{\bf{Proof of Theorem \ref{thm:VZW}}} 
For any non-negative integers $m,n$, from Corollary \ref{thm:main theorem2-2} we have
\begin{align}
\hat{A}^{m}A^{n}
   &=
   \begin{pmatrix}
   d_{m} & -q^{-m}b_{m} \\
   -q^{m}c_{m} & a_{m}
   \end{pmatrix}
   \begin{pmatrix}
   a_{n} & b_{n} \\
   c_{n} & d_{n}
   \end{pmatrix} \nonumber \\
\label{eq:VZW1}
   &=
   \begin{pmatrix}
   d_{m}a_{n}-q^{-m}b_{m}c_{n} & d_{m}b_{n}-q^{-m}b_{m}d_{n} \\
   -q^{m}c_{m}a_{n}+a_{m}c_{n} & -q^{m}c_{m}b_{n}+a_{m}d_{n}
   \end{pmatrix}
\end{align}
and 
\begin{align}
A^{n}\hat{A}^{m}
   &=
   \begin{pmatrix}
   a_{n} & b_{n} \\
   c_{n} & d_{n}
   \end{pmatrix}
   \begin{pmatrix}
   d_{m} & -q^{-m}b_{m} \\
   -q^{m}c_{m} & a_{m}
   \end{pmatrix} \nonumber \\
\label{eq:VZW2}
   &=
   \begin{pmatrix}
   a_{n}d_{m}-q^{m}b_{n}c_{m} & -q^{-m}a_{n}b_{m}+b_{n}a_{m} \\
   c_{n}d_{m}-q^{m}d_{n}c_{m} & -q^{-m}c_{n}b_{m}+d_{n}a_{m}
   \end{pmatrix}.
\end{align}
On the other hand, by applying $A\hat{A}=\hat{A}A=\delta E_{2}$ we obtain 
\begin{align}
\hat{A}^{m}A^{n}
   &=
   A^{n}\hat{A}^{m} \nonumber \\
\label{eq:VZW3}
   &=
\begin{cases}
\delta ^{m}A^{n-m}=
   \begin{pmatrix}
   \delta ^{m}a_{n-m} & \delta ^{m}b_{n-m} \\
   \delta ^{m}c_{n-m} & \delta ^{m}d_{n-m}
   \end{pmatrix} & (m< n) \\
\delta ^{n}\hat{A}^{m-n}=
   \begin{pmatrix}
   \delta ^{n}d_{m-n} & -q^{n-m}\delta ^{n}b_{m-n} \\
   -q^{m-n}\delta ^{n}c_{m-n} & \delta ^{n}a_{m-n}
   \end{pmatrix}
 & (m\geq  n)
\end{cases}.
\end{align}
By comparing the entries of (\ref{eq:VZW1}), (\ref{eq:VZW2}) and (\ref{eq:VZW3}) we have (\ref{eq:quantum rel1}), (\ref{eq:quantum rel2}), (\ref{eq:quantum rel3}) and (\ref{eq:quantum rel4}).

Finally, the relation (\ref{eq:quantum rel5}) follows from the explicit formulas (\ref{eq:explicit b}) and (\ref{eq:explicit c}):
\begin{align*}
b_{n}c_{m}
   &=
   q^{\frac{n-1}{2}}
   f_{n-1}(\tau^{\prime} )
   b
   q^{-\frac{m-1}{2}}c
   f_{m-1}(\tau ) \\
   &=
   q^{n-m}
   q^{-\frac{n-1}{2}}f_{n-1}(\tau^{\prime} )c
   q^{\frac{m-1}{2}}bf_{m-1}(\tau ) \\
   &=
   q^{n-m}c_{n}b_{m}.
\end{align*}
\qed

If we set $m=n$ in (\ref{eq:quantum rel1}), (\ref{eq:quantum rel2}), (\ref{eq:quantum rel3}), (\ref{eq:quantum rel4}) and (\ref{eq:quantum rel5}), then we obtain an interesting Corollary which means that $A^{n}$ is a $2\times 2$ $q^{n}$-quantum matrix.
\begin{cor}
\label{thm:cor q-rel}
For any non-negative integer $n$, we have
\begin{align}
a_{n}b_{n}&=q^{n}b_{n}a_{n}, \quad 
a_{n}c_{n}=q^{n}c_{n}a_{n}, \quad 
a_{n}d_{n}-d_{n}a_{n}=(q^{n}-q^{-n})b_{n}c_{n}, \nonumber \\
b_{n}c_{n}&=c_{n}b_{n}, \quad 
b_{n}d_{n}=q^{n}d_{n}b_{n}, \quad 
c_{n}d_{n}=q^{n}d_{n}c_{n} \tag{$R_{q^{n}}$}
\end{align}
and 
\begin{align}
a_{n}d_{n}-q^{n}b_{n}c_{n}
   =
   d_{n}a_{n}-q^{-n}b_{n}c_{n}
   =
   \delta ^{n}.
\end{align}
i.e. 
$$
(\text{the quantum determinant of $A^{n}$})=(\text{the quantum determinant of $A$})^{n}
$$
\end{cor}
Originally, Theorem \ref{thm:VZW} was proved by Vokos-Zumino-Wess \cite{VZW} and its proof was a brute force approach using double induction on $m$ and $n$. 
Later, Corrigan-Fairlie-Fletcher-Sasaki \cite{CFFS} and Umeda-Wakayama \cite{UW} gave some simple proofs of Corollary \ref{thm:cor q-rel} which is the case of $m=n$ in Theorem \ref{thm:VZW} independently. 
Our proof of Theorem \ref{thm:VZW} is different from any of them. 

In particular, it is desirable to extend Theorem \ref{thm:main theorem1} and Corollary \ref{thm:main theorem2} to $n\times n$ quantum matrices.

\section*{Acknowledgement}
We would like to thank Professor T\^{o}ru Umeda (Osaka City University) and Professor Masato Wakayama (Nippon Telegraph and Telephone Corporation, Institute for Fundamental Mathematics) for their valuable comments for references and various techniques of quantum calculus. 
This work was supported by Grant-in-Aid for Young Scientists (Number 21K13808).


\end{document}